\documentclass{article}
\usepackage{graphicx} 
\usepackage{amsmath}
\usepackage{amssymb}
\usepackage{indentfirst}
\usepackage{cite}
\usepackage{xcolor}
\usepackage[utf8]{inputenc}
\usepackage{subfigure}
\usepackage{float}
\usepackage{makeidx}  
\usepackage{hyperref}
\usepackage{amsthm}

\definecolor{akcolor}{rgb}{0.65, 0.15, 0.6}

\renewcommand{\qedsymbol}{$\blacksquare$} 

\title{Time-parameterized Optimal Transport}
\author{Kaiwen Shi}
\date{Fall 2024}

\begin{document}

\maketitle

\begin{abstract}
Optimal transport has gained significant attention in recent years due to its effectiveness in deep learning and computer vision. Its descendant metric, the Wasserstein distance, has been particularly successful in measuring distribution dissimilarities. While extensive research has focused on optimal transport and its regularized variants—such as entropy, sparsity, and capacity constraints—the role of time has been largely overlooked. However, time is a critical factor in real-world transport problems.

In this work, we introduce a time-parameterized formulation of the optimal transport problem, incorporating a time variable $t$ to represent sequential steps and enforcing specific constraints at each step. We propose a systematic method to solve a special subproblem and develop a heuristic search algorithm that achieves near-optimal solutions while significantly reducing computational time.
\end{abstract}

\section{Background}
The original optimal transport problem can be dated to the work of Gaspard Monge in 1781\cite{monge1781}. It was motivated by the transportation of masses between sources and sinks. Let $\delta$ denotes the Dirac delta function, and $\delta_{x}$ the  Dirac delta function whose mass concentrates on point $x$. Then we can let

\begin{equation} \label{source and sink}
    \mu = \sum_{i=1}^n a_i \delta_{x_i}, \quad \nu = \sum_{j=1}^m b_j \delta_{y_j}
\end{equation}
denote source and sink, respectively. When $\mu$ and $\nu$ have equal mass 1, that is, $\sum_{i=1}^n a_i=1=\sum_{j=1}^m b_j$,
the Monge formulation aims to find $T: \{x_1, x_2, \dots, x_n\} \to \{y_1, y_2, \dots, y_m\}$ such that $\forall j \in \{1, 2, \dots, m\}, b_j = \sum_{i: T(x_i) = y_j} a_i$.

To write a more elegant formulation of the Monge problem, we will need to define what a pushforward is. For $T$ a map from space $\mathcal{X}$ to space $\mathcal{Y}$, we define a corresponding pushforward operator  $T_\# :\mathcal{M}(\mathcal{X}) \to \mathcal{M}(\mathcal{Y})$ by $T_{\#} \mu (B) = \mu \left( \{ x: T(x) \in B\}\right), \forall B \subseteq Y$ a measurable set, where $\mathcal{M}(\mathcal{X})$ stands for the space of measures on $\mathcal{X}$. For discrete measures\footnote{For our purpose and interests, we will only consider discrete measures.}, the push-forward operation consists simply of moving the positions of all the points in the support of the measure. 
\begin{figure}[H]
    \centering
    \includegraphics[width=1\linewidth]{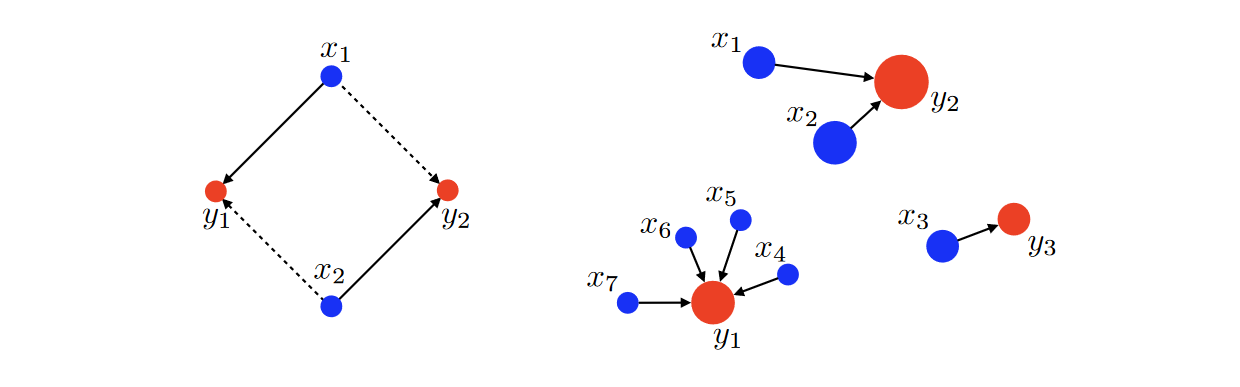}
    \caption{Visualizaion of pushforwards moving the masses at $x_i$ to position $y_j$\cite{peyre2019computational}.}
    \label{fig:Monge}
\end{figure}

Then, with this definition, the Monge problem is:
\begin{equation}
    \min_T \{\sum_i c(x_i, T(x_i)): T_\# \mu = \nu\}
\end{equation}
where $c(x,y)$ denotes the cost to transport mass from point $x$ to point $y$. In Euclidean spaces, it could be for example the Euclidean distance.

During WWII, the Soviet mathematician Leonid Kantorovich worked extensively on this problem and made huge contributions to this field \cite{kantorovich1942mass}. Part of his work included a reformulation of the transportation problem. Monge's problem might not have a solution because all the mass in one point from $\mathcal{X}$ will necessarily be transported to another point in $\mathcal{Y}$, with no room for moving the mass to other points. Kantorovich's relaxation of the problem gave a new perspective for splitting the mass: it allows mass from $x_i$ to be transported to more than one location.

Consider the discrete measures $\mu$ and $\nu$ as before. As they are discrete, we define vectors $\mathbf{a}$ and $\mathbf{b}$ with $\mathbf{a}_i = a_i$ and $\mathbf{b}_j = b_j$, essentially the coefficients of Dirac functions in $\mu$ and $\nu$ \eqref{source and sink}, to represent the measures $\mu$ and $\nu$. Then a transport plan $\mathbf{P}$ can be defined as

\begin{equation}
    \mathbf{P} \in \mathbb{R}_+^{n\times m} : \quad \mathbf{P}\mathbf{1}_m = \mathbf{a}, \quad \mathbf{P}^T\mathbf{1}_n = \mathbf{b}
\end{equation}
where $\mathbb{R}_+^{n \times m}$ denotes the space of non-negative real-valued matrices of size $n \times m$, and $\mathbf{1}_m$ denotes a column vector of size $m \times 1$ with all entries equal to 1, i.e., $\mathbf{1}_m = [1, 1, \dots, 1]^T$, and similarly, $\mathbf{1}_n$ denotes a column vector of size $n \times 1$ with all entries equal to 1.

Let $U(\mathbf{a},\mathbf{b})$ denote the set of all $\mathbf{P}$ satisfying the above condition. The Kantorovich's problem then is:

\begin{equation}
    L_c(\mathbf{a}, \mathbf{b})\quad:= \min_{\mathbf{P} \in U(\mathbf{a},\mathbf{b})} \langle \mathbf{C}, \mathbf{P}\rangle  \quad:= \sum_{i.j} \mathbf{C}_{i,k} \mathbf{P}_{i,j} \label{kantorovich}
\end{equation}
where the entries $\mathbf{C}_{i,j}$ represent the cost of transporting a unit mass from the $i$-th element of $\mu$ to the $j$-th element of $\nu$, $\mathbf{P}_{i,j}$ denotes the amount of mass transported from $x_i$ to $y_j$, and $\langle \mathbf{C}, \mathbf{P} \rangle$ represents the Frobenius inner product between the cost matrix $\mathbf{C}$ and the transport plan $\mathbf{P}$. 

Kantorovich later proved the existence of a solution \cite{kantorovich1942mass} and developed linear programming for finding the solution\footnote{ For curious readers, we include the reference to simplex algorithm\cite{simplex} and interior point method\cite{nemirovsky2004interior} for solving linear programming.}.In recognition of his contributions, Kantorovich was awarded the Nobel Memorial Prize in Economic Sciences in 1975, sharing the honor with Tjalling Koopmans. His methodologies have since become fundamental in optimization theory, influencing fields as diverse as mathematics, computer science, and operations research.

Kantorovich's work ultimately led to the formalization of the Wasserstein distance, the name of which honored the work of Russian mathematician Leonid Nisonovich Vaserstein, who in 1969 extended Kantorovich's ideas by defining a metric on the space of probability measures\cite{vaserstein1969markov}.

Also known as the Earth Mover's Distance (EMD), the Wasserstein distance is a metric that quantifies the cost of transporting one probability distribution to another. Formally, consider two probability measures $\mu$ and $\nu$ defined on a metric space $(\mathcal{X}, d)$, where $d$ is the distance function on $\mathcal{X}$. The $p$-Wasserstein distance, denoted as $W_p(\mu, \nu)$, is defined as:

\begin{equation}
    W_p(\mu, \nu) = \left( \inf_{\pi \in \Pi(\mu, \nu)} \int_{\mathcal{X} \times \mathcal{X}} d(x, y)^p \, d\pi(x, y) \right)^{1/p},
    \label{eq:wasserstein_distance}
\end{equation}
where $\Pi(\mu, \nu)$ is the set of all couplings of $\mu$ and $\nu$. A coupling $\pi \in \Pi(\mu, \nu)$ is a joint probability distribution over $\mathcal{X} \times \mathcal{X}$ with marginals $\mu$ and $\nu$.

The Wasserstein distance has become a cornerstone of optimal transport theory, with applications spanning machine learning, computer vision, and probability theory\cite{peyre2019computational,NIPS2015_a9eb8122,kolouri2019generalized}. Its intuitive interpretation as the minimal cost to "move mass" between distributions makes it particularly appealing in practical scenarios.

\section{Introduction}
Recent advancements have expanded the utility of Optimal Transport through techniques like entropy-regularized smoothing for efficient computation\cite{cuturi2013sinkhorn}, sparsity constraints for enhanced interpretability\cite{pmlr-v84-blondel18a}, and capacity constraints for realistic resource distribution\cite{korman2012optimaltransportationcapacityconstraints}. These developments highlight the versatility and power of OT in solving large-scale and high-dimensional problems, while ongoing research continues to explore new dimensions, such as temporal dynamics and their impact on transport processes.

Consider a supply chain problem where goods or supplies are transported by a fleet of vehicles. Transportation inherently takes time, and in many real-world scenarios, it is not feasible to complete the entire task of transporting goods in a single day. Instead, the transport plan must be adjusted dynamically across multiple days, reflecting changing constraints and demands. Moreover, demand patterns and supply availability can vary significantly over time, further complicating the planning process. This makes it critical to design transport plans that are not only optimal for individual time steps but also account for the cumulative effect of decisions over a longer horizon. By incorporating time into the framework of optimal transport, it becomes possible to better model real-world logistics and supply chain problems, where the efficient use of time and resources is paramount. Thus, a time-parameterized approach to optimal transport is both practically significant and a natural extension of classical formulations.

With this motivation, this paper aims to formulate a class of time-parameterized optimal transport problems that are closely aligned with real-world applications. Specifically, we incorporate two critical constraints into the problem: capacity constraints and sparsity constraints. Capacity constraints account for the total transportation capability available on a given day, considering the collective limitations of all vehicles. Sparsity constraints, on the other hand, limit the number of vehicles that can be deployed each day, reflecting operational or logistical restrictions.

We begin by formulating the time-parameterized optimal transport problem and then focus on two related subproblems, analyzing each in detail. For each subproblem, we propose a dedicated solution method. Finally, we discuss how these solutions can be systematically combined to address the overarching time-parameterized optimal transport problem under constraints of capacity and sparsity.

\section{Notation}
\begin{itemize}
    \item $\mu$, $\nu$ = (discrete probability) measures. Normally they assume the form of a linear combination of Dirac functions. (i.e. $\mu = \sum_1^n a_i\delta_{x_i}$ and $\nu = \sum_1^m a_j\delta_{y_j}$)
    \item $\mathbf{a}$, $\mathbf{b}$ = vector representation of $\mu$ and $\nu$. For $\mu = \sum_1^n a_i\delta_{x_i}$, $\mathbf{a} = (a_1, a_2, \dots, a_n)^T$. Similar for $\mathbf{b}$.
    \item $N$ = the total number of time steps.
    \item $i$ = index for a time step.
    \item $\gamma^i$ = transport plan at time step $i$ with $\gamma^i \in \mathbb{R}_+^{n\times m}$.
    \item $\mathbf{C}^i$ = cost matrix at time step $i$. $\mathbf{C}^i_{j,k}$ denotes the cost to transport a unit from $x_j$ to $y_k$ at time $i$.
    \item $\mathbf{M}^i$ = capacity constraint matrix at time $i$. A transport plan $\gamma^i$ should satisfy $0 \leq \gamma^i_{j,k} \leq \mathbf{M}^i_{j,k}$.
    \item $\mathbf{1}$ = column vector with all entries being 1. If a subscript exists, it denotes the length of the vector.
\end{itemize}

\section{Time Parameterized Optimal Transport}
We formulate in this section a time-parameterized optimal transport problem with two key constraints: capacity, which limits total daily transport capability, and sparsity, which restricts the number of deployable vehicles per day. These constraints reflect real-world operational challenges and extend classical formulations to address practical applications.

\begin{equation*}
\begin{aligned}
\underset{\gamma^1, \dots, \gamma^N}{\arg\min} \quad & \sum_{i=1}^N \langle \mathbf{C}^i, \gamma^i \rangle = \sum_{i=1}^N \sum_{j,k}\mathbf{C}^i_{j,k} \gamma^i_{j,k} \\
\text{subject to} \quad & 0 \leq \gamma^i_{j,k} \leq \mathbf{M}^i_{j,k}, \\
                        & \|\gamma_{j,\cdot}^i\|_0 \leq s_{j}^i, \\
                        & \mathbf{F} = \sum_i \gamma^i, \\
                        & \mathbf{F}\mathbf{1} = \mathbf{a}, \quad \mathbf{F}^T\mathbf{1} = \mathbf{b}
\end{aligned}
\label{Merged}
\end{equation*}
Here:
\begin{itemize}
    \item $\gamma^i$ represents the transport plan at time step $i$.
    \item $\mathbf{C}^i$ denotes the cost matrix for time step $i$.
    \item $\mathbf{M}^i$ is the capacity matrix for time step $i$, imposing the upper bounds on transport capabilities.
    \item $s^i_{j}$ represents the sparsity constraint, limiting the number of vehicles from source $j$.
    \item $\mathbf{F}$ is the cumulative transport plan across all time steps. It is computed to ensure the match of the marginals.
    \item $\mathbf{a}$ and $\mathbf{b}$ are the source and sink probability measures, respectively.
\end{itemize}

This formulation extends classical optimal transport by integrating temporal dynamics and practical constraints, making it highly relevant for real-world logistics. However, the broad scope of this formulation introduces significant challenges: it is inherently nonconvex and computationally complex, making it difficult to analyze and solve directly. To address these issues, we decompose the problem into two relatively more manageable subproblems. The first focuses on \textbf{time-parameterized optimal transport with capacity constraints}, while the second examines \textbf{sparsity-constrained optimal transport}\footnote{For the second subproblem, the time variable is excluded, as time steps are meaningful only when daily capacity is limited. Otherwise, the entire transportation can be completed in a single day.}.

\textbf{Note:} Curious readers may notice that the formulation assumes the source and sink are probability measures, which is not always the case in real-world scenarios. In practice, the total mass of the measures can exceed 1. However, this issue can be resolved through normalization by a constant. The following lemma establishes that if the measures are scaled by a constant factor, the optimal transport plan is similarly scaled. Consequently, to recover the original optimal transport plan before normalization, it suffices to reapply the constant multiplier.

\underline{\textbf{Lemma:}} A scaled fraction of an optimal transport plan remains optimal for correspondingly scaled marginals. Specifically, let $P$ be an optimal transport plan between $\mathbf{a}$ and $\mathbf{b}$, and let $c > 0$ be a positive scalar with $\frac{1}{c} < 1$. Then, $\frac{P}{c}$ is an optimal transport plan between $\frac{\mathbf{a}}{c}$ and $\frac{\mathbf{b}}{c}$.

\underline{\textbf{Proof:}} We prove this lemma by contradiction. Assume $P$ is an optimal transport plan between $\mathbf{a}$ and $\mathbf{b}$, but $\frac{P}{c}$ is not optimal between $\frac{\mathbf{a}}{c}$ and $\frac{\mathbf{b}}{c}$. Then, there must exist another transport plan $P^*$ that is optimal between $\frac{\mathbf{a}}{c}$ and $\frac{\mathbf{b}}{c}$. Multiplying $P^*$ by $c$ yields a transport plan $cP^*$ between $\mathbf{a}$ and $\mathbf{b}$. However, this plan $cP^*$ has a lower cost than $P$ since $P^*$ has a lower cost than $\frac{P}{c}$, which contradicts the assumption that $P$ is optimal between $\mathbf{a}$ and $\mathbf{b}$. Hence, the lemma holds. \qedsymbol

\section{Time-parameterized Optimal Transport with Capacity Constraints}
\label{TOTCC}
\subsection{Problem Statement}
Given two discrete probability measures $\mu$ and $\nu$, we aim to find a sequence of transport matrices $\gamma^i$, $i \in \{1, 2, \ldots, N\}$, such that for each $i$, $\gamma^i$ represents a transport plan between some measures $\mu^i$ and $\nu^i$, the marginals of each time step. These measures satisfy the constraints $\int d\mu^i \leq \int d\mu$ and $\int d\nu^i \leq \int d\nu$, with $\mu$ and $\nu$ serving as the marginals of the summative transport plan $(\gamma^1 + \ldots + \gamma^N)$. Additionally, the transport plans $\gamma^i$ must adhere to specific constraints, such as capacity constraints.

For discrete measures $\mu = \sum_{i=1}^n a_i\delta_{x_i}$ and $\nu = \sum_{j=1}^m b_j\delta_{y_j}$, where $\mathbf{a} = (a_1, \ldots, a_n)^T$ and $\mathbf{b} = (b_1, \ldots, b_m)^T$ represent the vectors of the weights, the time-parameterized optimal transport problem with capacity constraints can be formulated as the following optimization problem:

\begin{equation*}
\begin{aligned}
\underset{\gamma^1, \dots, \gamma^N}{\arg\min} \quad & \sum_{i=1}^N \langle \mathbf{C}^i, \gamma^i \rangle = \sum_{i=1}^N \sum_{j,k} \mathbf{C}^i_{j,k} \gamma^i_{j,k} \\
\text{subject to} \quad & 0 \leq \gamma^i_{j,k} \leq \mathbf{M}^i_{j,k}, \\
                        & \mathbf{F} = \sum_i \gamma^i, \\
                        & \mathbf{F}\mathbf{1} = \mathbf{a}, \quad \mathbf{F}^T\mathbf{1} = \mathbf{b}.
\end{aligned}
\end{equation*}
where $\mathbf{M}^i$ is the capacity matrix for time step $i$, imposing the upper bounds on transport capabilities.

Notice that for all $j$, there exists at least an $i$ with $\sum_k M^i_{j,k} > 0$. That is to say, in a more practical term, for a specific source, it cannot be that every $M^i_{j,k}$ has all zero entries, which in paraphrase means the capacity to transport from this source is zero. 

\underline{\textbf{Lemma:}}  
For the time-parameterized optimal transport problem, assume $\mathbf{M}^i_{j,k}$ represents the capacity constraints for transport plans $\gamma^i$. If for every source index $j$, there exists at least one time step $i$ such that $\sum_k \mathbf{M}^i_{j,k} > 0$, then a feasible solution to the problem can exist. Conversely, if for any $j$, $\sum_k \mathbf{M}^i_{j,k} = 0$ for all $i$, it implies that no transport capacity is available from source $j$, and thus no feasible solution exists. 

\underline{\textbf{Proof:}} The proof is rather trivial. If for any $j$, $\sum_k \mathbf{M}^i_{j,k} = 0$ for all $i$, then $(\mathbf{\sum_i \gamma^i}\mathbf{1})_j \leq (\mathbf{\sum_i \mathbf{M}^i}\mathbf{1})_j = (\sum_i\sum_k \mathbf{M}^i_{j,k})_j = 0 < \mathbf{a}_j$, which does not satisfy the constraint $ \mathbf{\sum_i \gamma^i}\mathbf{1} = \mathbf{F}\mathbf{1} =\mathbf{a}$. \qedsymbol

A similar proof can be given to show that for all $k$, there exists at least an $i$ with $\sum_j \mathbf{M}^i_{j,k} > 0$. Again, this means for some days, the transport capacity to sink $k$ has to be non-zero. Otherwise, the marginal constraint would not be met for $\mathbf{b}$.

Note that the above problem is written for a fixed number $N$ of time steps. However, one can also let $N$ be free and consider a new problem where the total number of time steps used to partition the final plan $F$ is not determined \textit{a priori}. In this new scenario, we have to minimize over all the possible finite partitions $\{\gamma^i\}$ of $\mathbf{F}=\sum_i\gamma^i$. We will not include any further discussion on this new problem since it is out of the scope of this paper.

\subsection{A Special Case}
The above problem can be solved using linear programming. However, in a special case where the cost matrix and capacity constraints are constant across time steps, we can derive an equivalent formulation and reduce the computational cost.

The special case occurs when the cost matrix and capacity constraint matrix are the same for all time steps. That is, $\mathbf{C}^i = \mathbf{C}^j$ and $\mathbf{M}^i = \mathbf{M}^j$ for all $i,j \in \{1, \dots, N\}, i \neq j$. In this case, the original formulation simplifies to:
\begin{equation} 
\begin{aligned}
\underset{\gamma^1, \dots, \gamma^N}{\arg\min} \quad & \langle \mathbf{C}, \mathbf{F} \rangle =  \sum_{j,k} \mathbf{C}_{j,k} \mathbf{F}_{j,k} = \sum_{i=1}^N \sum_{j,k}\mathbf{C}_{j,k} \gamma^i_{j,k} \label{eq:special case}\\
\text{subject to} \quad & 0 \leq \gamma^i_{j,k} \leq \mathbf{M}_{j,k}, \qquad \text{(capacity constraint)}\\
                        & \mathbf{F} = \sum_i \gamma^i, \\
                        & \mathbf{F}\mathbf{1} = \mathbf{a}, \quad \mathbf{F}^T\mathbf{1} = \mathbf{b}
\end{aligned}
\end{equation}

or simply as
\begin{equation}\label{eq: problem C M const}
    \min_{\mathbf{F} \in U_M(\mathbf{a},\mathbf{b})} \langle \mathbf{C}, \mathbf{F}\rangle  
\end{equation}
where the set of feasible transportation plans for the above problem is
\begin{equation}
  U_M(\mathbf{a}, \mathbf{b}):=\left\{\mathbf{F}: \, \mathbf{F} = \sum_i \gamma^i\in U(\mathbf{a}, \mathbf{b}), \text{ with } 0 \leq \gamma^i_{j,k} \leq \mathbf{M}_{j,k}, \,
                      \right\}  \label{Feasible Set}
\end{equation}

\subsubsection{An Example of the Special Case and Solution}  
To help readers better understand the problem, we present an example. Suppose we have two mines and two warehouses, and the task is to complete the transport between them within two days. Mine 1 has 6 units of products, and Mine 2 has 8 units. We want Warehouse 1 to be stocked with 4 units, and Warehouse 2 with 10 units, by the end of day 2. Each transport incurs some cost and is subject to capacity constraints.

Let $\mu = 6\delta_{x_1} + 8\delta_{x_2}$ and $\nu = 4\delta_{y_1} + 10\delta_{y_2}$ be the measures, with the cost matrix $\mathbf{C} = \begin{pmatrix} 1 & 4 \\ 3 & 6 \end{pmatrix}$ and the capacity constraint matrix $\mathbf{M} = \begin{pmatrix} 1 & 2 \\ 2 & 4 \end{pmatrix}$ (i.e., $0 \leq \gamma^i_{j,k} \leq \mathbf{M}_{j,k}$).

We seek to solve the following optimization problem:
\begin{equation*}
\begin{aligned}
\underset{\gamma^1,\gamma^2}{\arg\min} \quad & \langle \mathbf{C}, \mathbf{F} \rangle =  \sum_{j,k} \mathbf{C}_{j,k} \mathbf{F}_{j,k} = \sum_{i=1}^2 \sum_{j,k} \mathbf{C}_{j,k} \gamma^i_{j,k} = \sum_{i=1}^2(\gamma^i_{1,1} + 4\gamma^i_{1,2} + 3\gamma^i_{2,1} + 6\gamma^i_{2,2}) \\
\text{subject to} \quad & 0 \leq \gamma^i_{j,k} \leq \mathbf{M}_{j,k}, \qquad \text{(capacity constraint)} \\
                        & \mathbf{F} = \sum_i \gamma^i, \\
                        & \mathbf{F}\mathbf{1} = \mathbf{a} = (6,8), \quad \mathbf{F}^T\mathbf{1} = \mathbf{b} = (4,10)
\end{aligned}
\end{equation*}

\underline{Solution:}  
$\gamma^1 = \begin{pmatrix} 1 & 2 \\ 1 & 4 \end{pmatrix}$, $\gamma^2 = \begin{pmatrix} 1 & 2 \\ 1 & 2 \end{pmatrix}$, and $\mathbf{F} = \begin{pmatrix} 2 & 4 \\ 2 & 6 \end{pmatrix}$. (Note that the solution may not be unique.)

\subsubsection{Observation:}  
\label{subsubsec:obs}  
Let $\mathbf{F}^\dagger$ be the optimal solution for the problem \eqref{eq:special case}, and let $\mathbf{F}^*$ be the optimal transport plan for the Kantorovich problem \eqref{kantorovich}:

\begin{equation*}
\begin{aligned}
\underset{\mathbf{F} \in U(\mathbf{a,b})}{\arg\min} \quad & \langle \mathbf{C}, \mathbf{F}\rangle  = \sum_{j,k} \mathbf{C}_{j,k} \mathbf{F}_{j,k} \\
\text{subject to} \quad & \mathbf{F}\mathbf{1} = \mathbf{a}, \quad \mathbf{F}^T\mathbf{1} = \mathbf{b}
\end{aligned}
\end{equation*}

Then, the following inequality holds:

\begin{equation} \label{eq: less or equal}
    \langle \mathbf{C}, \mathbf{F}^*\rangle \leq \langle \mathbf{C}, \mathbf{F}^\dagger \rangle
\end{equation}

This observation stems from the fact that the problem \eqref{eq:special case} is a restricted version of the original Kantorovich problem, operating on a more limited feasible set, i.e., $U_M(\mathbf{a},\mathbf{b}) \subseteq U(\mathbf{a},\mathbf{b})$, where $U(\mathbf{a},\mathbf{b})$ is defined as the set of all solutions as in Kantorovch \eqref{kantorovich}. As a result, the optimized cost in the restricted problem cannot beat that of the Kantorovich solution.

\subsubsection{Existence of Solution(s):}  
We now consider the existence of solutions to the problem \eqref{eq:special case}. Due to the following proposition, we will focus on the case where $N$, the total number of time steps, is fixed.

\underline{\textbf{Proposition:}}  
If $N$ is free and $\mathbf{M}$ has all non-zero entries, then $\mathbf{F}^\dagger = \mathbf{F}^*$.

\underline{\textbf{Proof:}}  
Since $N$ is free, we can choose an arbitrarily large $N$. Let $N = \left\lceil \max_{j,k} \left( \frac{\mathbf{F}^*_{j,k}}{\mathbf{M}_{j,k}} \right) \right\rceil$. Then, define $\gamma^i = \frac{\mathbf{F}^*}{N}$ for all $i$. This can be done without violating any constraints, since $N \geq \frac{\mathbf{F}^*_{j,k}}{\mathbf{M}_{j,k}}$ ensures that $\gamma^i_{j,k} = \frac{\mathbf{F}^*_{j,k}}{N} \leq \mathbf{M}_{j,k}$. Notice that $\sum_{i=1}^N \gamma^i = \mathbf{F}^*$, and the solution satisfies the constraints of the problem \eqref{eq:special case}. Therefore, $\mathbf{F}^* = \mathbf{F}^\dagger$ due to observation \eqref{eq: less or equal}. \qedsymbol

With $N$ fixed, we have found a necessary and sufficient condition for the existence of solutions to the problem \eqref{eq:special case}. The proof follows similar logic found in the literature on Capacity Constrained Optimal Transport, such as \cite{korman2012optimaltransportationcapacityconstraints}.

\label{existence}
\underline{\textbf{Proposition:}} (Existence of Solutions) If the feasible set of transportation plans \eqref{Feasible Set} is non-empty, then a solution(minimizer) exists. 

\underline{\textbf{Proof:}} We prove by showing that the domain of the function we wish to minimize defined in \eqref{eq:special case} is compact. First, we note that for a sequence of transport plans $\gamma^1, \dots, \gamma^N$, one can view them as a sequence of $\mathbb{R}^{n \times m}_+$ matrices, or equivalently, a vector in $\mathbb{R}^{N \cdot n \cdot m}_+$ space.

Now we focus on the compactness of the domain of the function. The linear inequalities then define a high-dimensional cube that is both closed and bounded. (That is, $[0,\mathbf{M}_{1,1}] \times [0,\mathbf{M}_{1,2}] \times \dots \times [0,\mathbf{M}_{m,n}]$) Also, the linear inequalities essentially encompass an intersection of affine hyperplanes, thus closed. Since a closed subset of a closed and bounded set is closed and bounded, we see the constrained domain is indeed closed and bounded. Then by the Heine-Borel Theorem, the domain is compact.

Since by the assumption we have a non-empty solution set (a very common assumption in Capacity Constrained Optimal Transport), a minimizer thus exists. \qedsymbol

\textbf{Note: }One way to ensure that the domain is not empty is to enforce one possible candidate plan. For instance, if we enforce this condition that $\mathbf{M}_{j,k} \geq \mathbf{P}_{j,k}, \forall j,k$, where $\mathbf{P} = \mathbf{a}^T \mathbf{b}$, the domain must contain the sequence of fractioned product measure $\gamma^i = \frac{\mathbf{a}^T \mathbf{b}}{N}, \forall i$.

\subsection{A Faster Algorithm to Compute Solution}

In this section, we propose an alternative yet faster algorithm to solve the problem \eqref{eq:special case}. 

\subsubsection{Observation:}
\label{ReformulationOfProb}
It is very helpful to note that solving the problem with free $\gamma^i$'s is equivalent to solving the problem on a more restricted feasible set \eqref{Feasible Set} with $\gamma^i$'s equal. Indeed, we will show that an equivalent formulation for problem in \eqref{eq:special case}:

\begin{equation*}
\begin{aligned}
\underset{\widehat{\gamma}}{\arg\min} \quad & \langle \mathbf{C}, \mathbf{P} \rangle =  \sum_{j,k} \mathbf{C}_{j,k} \mathbf{P}_{j,k} 
= N \sum_{j,k}\mathbf{C}_{j,k} \widehat{\gamma}_{j,k} \\
\text{subject to} \quad 
& 0 \leq \widehat{\gamma}_{j,k} \leq \mathbf{M}_{j,k}, \\
& \mathbf{P} = N \cdot \widehat{\gamma}, \\
& \mathbf{P}\mathbf{1} = \mathbf{a}, \quad \mathbf{P}^T\mathbf{1} = \mathbf{b}.
\end{aligned}
\end{equation*}

In other words, without loss of generality, we can assume \(\gamma^i = \gamma^j, \forall i \neq j\). 

The interpretation of this reads as follows: with cost and capacity constraint matrices being unchanged, the transport plans across all time steps should not be varied as they are weighted equally. Precisely, we provide the following result.

\underline{\textbf{Proposition:}}
Consider 
\begin{equation}
    U_M^{Unif}(\mathbf{a},\mathbf{b}):=\left\{\mathbf{P}: \, \mathbf{P}=N \cdot \widehat{\gamma}\in U(\mathbf{a},\mathbf{b}),  \text{ with } 0\leq \widehat{\gamma}_{j,k}\leq \mathbf{M}_{j,k}\right\}
\end{equation}

Then

\begin{equation}
\inf_{\mathbf{F}\in\ U_M(\mathbf{a},\mathbf{b})}\langle \mathbf{C},\mathbf{F}\rangle =\inf_{\mathbf{P}\in\ U_M^{Unif}(\mathbf{a},\mathbf{b})}\langle \mathbf{C},\mathbf{P}\rangle 
\end{equation}

\underline{\textbf{Proof:}} To show this, we will show
\[
\underset{\mathbf{F}\in U_M(\textbf{a},\textbf{b})}{\inf} \langle \mathbf{C}, \mathbf{F} \rangle \leq \underset{\mathbf{P}\in U_M^{Unif}(\textbf{a},\textbf{b})}{\inf} \langle \mathbf{C}, \mathbf{P} \rangle \quad \text{and} \quad \underset{\mathbf{F}\in U_M(\textbf{a},\textbf{b})}{\inf} \langle \mathbf{C}, \mathbf{F} \rangle \geq \underset{\mathbf{P}\in U_M^{Unif}(\textbf{a},\textbf{b})}{\inf} \langle \mathbf{C}, \mathbf{P} \rangle
\]

The left inequality is obvious. The logic follows what has been written in observation \ref{subsubsec:obs}:  as $U_M^{Unif}(\mathbf{a},\mathbf{b})$ is a more restricted set (that is, $U_M^{Unif}(\mathbf{a},\mathbf{b})\subseteq U_M(\mathbf{a},\mathbf{b})$), the optimal solution cannot get better.

Now we will prove the right inequality. Let $\mathbf{F} = \sum_i \gamma^i$ be an arbitrary optimal solution plan for \eqref{eq: problem C M const}, and let $\widehat{\gamma} = \frac{\sum_i \gamma^i}{N}$, and $ \mathbf{P}:= N \cdot \widehat{\gamma}$.

We will need to show that $\widehat{\gamma}$ constitutes a feasible solution first. Since $0 \leq \gamma^i_{j,k} \leq \mathbf{M}_{j,k}$, we have $0 \leq \sum_i \gamma^i_{j,k} = N \cdot \widehat{\gamma}_{j,k} \leq N\cdot \mathbf{M}_{j,k}$, thus $0 \leq \widehat{\gamma}_{j,k} \leq \mathbf{M}_{j,k}$. Therefore, the capacity constraint is satisfied.

Also, we show that $P$ is no worse than $F$ in terms of optimality. Indeed,
\[
\langle \mathbf{C}, \mathbf{F} \rangle = \langle \mathbf{C}, \sum_i \gamma^i \rangle = \langle \mathbf{C}, N \cdot \frac{\sum_i \gamma^i}{N} \rangle = \langle \mathbf{C}, N \cdot \widehat{\gamma} \rangle =  \langle \mathbf{C}, \mathbf{P} \rangle
\]

Thus, the equality of optimal cost holds. Indeed, we can reformulate the problem on a more restricted set and assume all $\gamma^i$'s are equal. \qedsymbol

\subsubsection{A Faster Algorithm}
Here, we introduce a faster algorithm to solve the problem \eqref{eq:special case} by leveraging the reformulation in observation \ref{ReformulationOfProb}. Due to the reformulation, the dimensionality of the solution space is significantly reduced. Previously, we established that the original solutions can be represented as vectors in $\mathbb{R}^{N \cdot n \cdot m}_+$, where $N$ is the number of time steps, and $n$ and $m$ are the dimensions of the problem. With the reformulation in observation \ref{ReformulationOfProb}, the solutions are now restricted to vectors in $\mathbb{R}^{n \cdot m}_+$. This dimensionality reduction leads to computational advantages.

In the current literature and practice \cite{BPPH11}, solving a linear programming problem typically has a time complexity of $O(n^3L)$, where $L = \log(n)$, and $n$ is the number of variables. Sometimes this complexity is simplified to $\Tilde{O}(n^3)$, omitting the logarithmic term. By applying our reformulation, the theoretical time complexity can be reduced by a factor of $O(N^3\log(N))$ under the assumption that $n = m$. In simpler terms, this means that reducing the size of the problem (by decreasing the number of variables) allows the algorithm to solve it more quickly.

Other previously used methods, such as the simplex algorithm and interior-point methods, have higher time complexities. The simplex method is known to be weakly polynomial \cite{dantzig1998linear}, while interior-point methods have a time complexity of approximately $O(n^{3.5}L)$ \cite{karmarkar1984new}. In these cases, our reformulated approach provides even greater time-saving benefits.

\subsubsection{Experiment}
To illustrate the improvement in computation efficiency, we conduct an experiment with Python. For the experiment, we use randomly generated $\mu \in \mathbb{R}^n_+$, $\nu \in \mathbb{R}^m_+$, and $\mathbf{C} \in \mathbb{R}^{n\times m}_+$. We also randomly generate $\mathbf{M} \in \mathbb{R}^{n\times m}_+$ such that condition \ref{existence} is met (i.e. a solution exists). Then we use different solvers to solve the linear programming problem and document the results. We take the average of 10 runs, but in cases where an individual run takes too little time, the results are taken as the average of 100 successful runs to mitigate unwanted software/hardware perturbations.

\subsubsection{Result}
The results of the experiment are as follows\footnote{The codes for this experiment can be found at this Notebook: \href{https://colab.research.google.com/drive/1ll-hGtKwlBl8SMSx7VIYhjG3tbDRAg9i?usp=sharing}{Colab Notebook link}. Note HiGHS stands for High-performance parallel linear optimization software. It deploys the power of parallel computing and problem sparsity to speed up computation.}.

\begin{table}[H]
\centering
\title{\textbf{Table 1: Experiment Results when n = m = 10 but N varies}}
\begin{tabular}{|c|c|c|c|c|}
\hline
\textbf{Method} & \textbf{Revised Simplex} & \textbf{Interior Point} & \textbf{HiGHS} \\ 
\hline
{N = 10}    &   3.04 &  3.95  &  0.0072 \\ 
\hline
\textbf{N = 10}  &  \textbf{0.052}  &  \textbf{0.045}  & \textbf{0.0047}   \\ 
\hline
{N = 50}    & 196.12   & 226.28   &  0.038  \\ 
\hline
\textbf{N = 50}   &  \textbf{0.054}  &  \textbf{0.056}  & \textbf{0.0047}   \\ 
\hline
{N = 100}    & 1417.64   &   1925.13 &  0.058\\ 
\hline
\textbf{N = 100}   &    \textbf{0.087} &  \textbf{0.064}  & \textbf{0.0074}   \\ 
\hline
\end{tabular}
\caption{The bold fonts in the rows denote the speed of the reformulated problem, compared to normal fonts which indicate the speed of the original formulation. Notice that the bold fonts in the same column should be similar.}
\end{table}

\section{Sparsity Constrained Optimal Transport}
The second subproblem we aim to address is the optimal transport problem with sparsity constraints. To motivate this again, consider a practical example: suppose you need to transport goods from several sources to multiple destinations. However, there is a limit on the number of transportation vehicles available at each source (e.g., three trucks per day per source). Consequently, the transport plan must ensure that only a limited number of routes, represented by non-zero entries in the transport plan, are utilized from each source. 

In recent studies \cite{blondel2018smoothsparseoptimaltransport, liu2023sparsityconstrainedoptimaltransport}, it has been demonstrated that while quadratic regularization can yield sparse transport plans, it does not directly control the exact number of non-zero entries. Methods have been developed to explicitly limit this number to impose stricter and more precise sparsity constraints.

The problem is formulated as follows:
\begin{equation*}
\begin{aligned}
\underset{\gamma}{\arg\min} \quad & \langle \mathbf{C}, \mathbf{\gamma} \rangle = \sum_{j,k} \mathbf{C}_{j,k} \gamma_{j,k}, \\[0.5em]
\text{subject to} \quad & \|\gamma_{j,\cdot}\|_0 \leq s_j, \quad \forall j \quad \text{(sparsity constraint)}, \\[0.5em]
                        & \mathbf{\gamma}\mathbf{1} = \mathbf{a}, \quad \mathbf{\gamma}^T\mathbf{1} = \mathbf{b},
\end{aligned}
\label{eq:sparsePrimal}
\end{equation*}
where:
\begin{itemize}
    \item $\mathbf{C} \in \mathbb{R}^{n \times m}_+$ is the cost matrix, 
    \item $\|\cdot\|_0$ is the $\ell_0$ norm (counting the number of non-zero elements), and 
    \item $s_j \in \mathbb{R}$ specifies the upper bound on the number of non-zero entries for each source $j$.
\end{itemize}

An attentive reader might notice that this formulation excludes capacity constraints, instead focusing exclusively on sparsity constraints. This simplification is intentional, as it isolates the sparsity aspect for clarity and precludes the necessity of the use of the time parameter $t$.

\subsection{Non-convexity Poses Difficulty}
One notable difficulty with this formulation lies in the use of the $\ell_0$ norm, which is non-convex and computationally challenging to optimize. A standard approach to make the problem with $\ell_0$ norm more tractable is to relax the $\ell_0$ norm to the $\ell_1$ norm. 

As explored in \cite{ramirez2013why}, the $\ell_1$ norm is a natural convex approximation to the $\ell_0$ norm. The authors offer a geometric perspective, defining a partial order $\leq$ on functions, represented by a set 
\[
B_{\leq} = \{c : c \leq (1, 0, \dots, 0)\} = \{c : f(c_1, \dots, c_n) \leq f(1, 0, \dots, 0)\}.
\]
They demonstrated that every convex set containing $B_0$ (associated with the $\ell_0$ norm) must also contain $B_1$ (associated with the $\ell_1$ norm). Hence, the $\ell_1$ norm is geometrically the closest convex relaxation of the $\ell_0$ norm.

\begin{figure}[H]
    \centering
    \subfigure[$B_0 = \bigcup_{i=1}^n \{(0, \dots, c_i, \dots, 0) : |c_i| \leq 1\}$]{
        \includegraphics[width=0.45\linewidth]{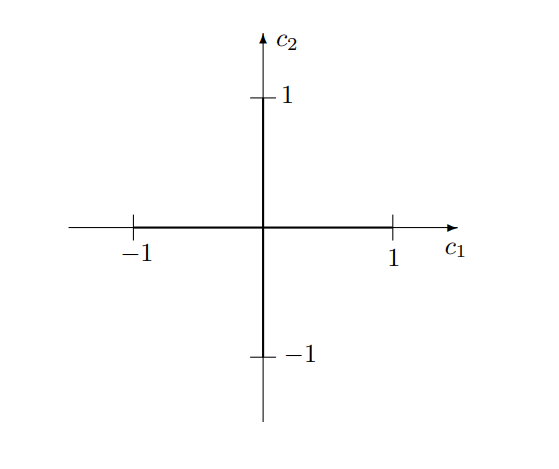}
        \label{fig:subfig1}
    }
    \hspace{0.05\linewidth}
    \subfigure[$B_1 = \{(c_1, \dots, c_n): \sum_{i=1}^n |c_i| = 1\}$]{
        \includegraphics[width=0.45\linewidth]{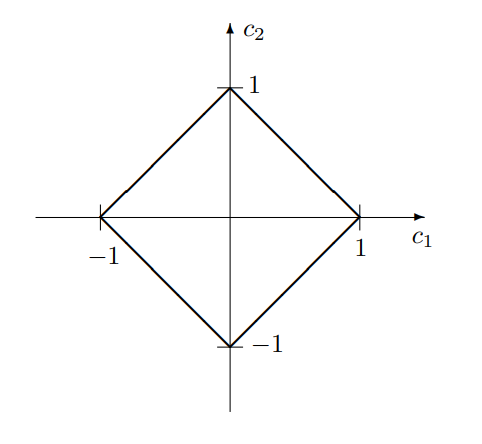}
        \label{fig:subfig2}
    }
    \caption{The sets $B_0$ and $B_1$ as defined in the paper. $B_1$ is the convex hull of $B_0$.}
    \label{fig:combined}
\end{figure}

Another perspective on why the $\ell_1$ norm promotes sparsity can be found in the LASSO paper \cite{tibshirani1996regression}. The authors reformulate the regression problem 
\[
\underset{\hat{\alpha},\hat{\beta}}{\arg\min} \quad \sum_{i=1}^n \left(y_i - \alpha - \sum_{j} \beta_j x_{ij}\right)^2, \quad \text{subject to } \sum_j |\beta_j| \leq t,
\]
as the unconstrained optimization problem:
\[
\underset{\hat{\alpha},\hat{\beta}}{\arg\min} \quad \sum_{i=1}^n \left(y_i - \alpha - \sum_{j} \beta_j x_{ij}\right)^2 + \lambda \sum_j |\beta_j|,
\]
where $\lambda$ controls the penalty strength.

The sparsity-inducing behavior of LASSO stems from its constraint region, which is a rotated square (a diamond in 2D). Unlike the spherical constraint of ridge regularization, the sharp corners of this square align with the coordinate axes, encouraging some coefficients $\beta_j$ to be exactly zero. This geometric property makes LASSO effective in promoting sparsity and selecting features in high-dimensional problems.

\begin{figure}[H]
    \centering
    \includegraphics[width=0.8\linewidth]{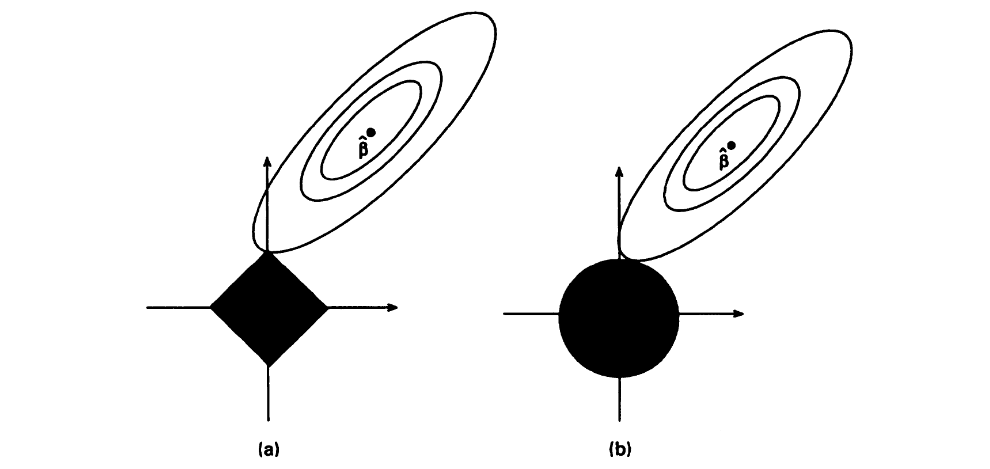}
    \caption{Left: LASSO. Right: Ridge. The contour is the level set of the cost function, and the shaded area is the constraint region. It is apparent that it is easier to have a tangent line of the level set going through a corner (a sparse vector) on the left than on the right.}
    \label{fig:enter-label}
\end{figure}

Building on these observations and empirical successes \cite{candes2006stable, candes2005decoding, donoho2006compressed, sanchez2011sparse, elad2010sparse}, we try to relax the $l_0$ norm to the $l_1$ norm\footnote{Note that this relaxation does not guarantee recovery of the original sparse solution.}. This gives the following primal formulation:

\[
\begin{aligned}
\underset{\gamma}{\arg\min} \quad & \langle \mathbf{C}, \mathbf{\gamma} \rangle = \sum_{j,k} \mathbf{C}_{j,k} \gamma_{j,k}, \\
\text{subject to} \quad & \|\gamma_{j,\cdot}\|_1 \leq s_j, \forall j \text{ (sparsity constraint)}, \\
                         & \gamma_{j,k} \geq 0, \forall j, k, \\
                         & \mathbf{\gamma} \mathbf{1} = \mathbf{a}, \quad \mathbf{\gamma}^T \mathbf{1} = \mathbf{b}.
\end{aligned}
\]

Although the relaxation makes the problem convex and solvable, it introduces issues. Due to the marginal constraints, $\|\gamma_{j,\cdot}\|_1$ is fixed at $\mathbf{a}_j$. If $s_j < \mathbf{a}_j$, the marginal constraints are infeasible, and if $s_j \geq \mathbf{a}_j$, the sparsity constraint becomes redundant. Consequently, the $l_1$ norm relaxation is not meaningful in this context, and relaxations with other norms fail to recover sparse solutions.

\subsection{A Proposed Heuristics}
The zero-norm formulation is a necessity, but solving the original problem \eqref{eq:sparsePrimal} is computationally expensive due to its combinatorial nature. Specifically, solving it requires selecting some entries to be zeros in the transport plan, iterating over all valid selections, and comparing them to find the optimal solution. If a systematic way to find the solution in one pass exists, eliminating the need to test all valid plans, it could significantly accelerate the process and reduce reliance on relaxations.

A heuristic approach involves defining an importance score for each entry in the transport plan and selecting non-zero entries based on this score. If an entry has a higher importance score, it is prioritized in the solution. For example, consider the sparse OT problem:

\[
\begin{aligned}
\underset{\gamma}{\min} \quad & \langle \mathbf{C}, \mathbf{\gamma} \rangle, \\
\text{subject to} \quad & \|\gamma^i_{j,\cdot}\|_0 \leq 2, \forall i, j, \\
                         & \mathbf{\gamma}\mathbf{1} = \mathbf{a}, \quad \mathbf{\gamma}^T\mathbf{1} = \mathbf{b},
\end{aligned}
\]
where $\mathbf{a}, \mathbf{b} \in \mathbb{R}^3$, $\mathbf{C} \in \mathbb{R}^{3 \times 3}$, and an importance score matrix:

\[
\mathbf{I} = 
\begin{bmatrix}
0.2 & 0.5 & 0.3 \\
0.1 & 0.6 & 0.3 \\
0.4 & 0.2 & 0.4
\end{bmatrix}.
\]

We know that $\gamma \in \mathbb{R}^{3 \times 3}$. A method that selects non-zero entries based on indices would first test:

\[
\begin{aligned}
\underset{\gamma}{\min} \quad & \langle \mathbf{C}, \mathbf{\gamma} \rangle, \\
\text{subject to} \quad & \gamma_{1,3} = \gamma_{2,3} = \gamma_{3,3} = 0, \\
                         & \mathbf{\gamma}\mathbf{1} = \mathbf{a}, \quad \mathbf{\gamma}^T\mathbf{1} = \mathbf{b},
\end{aligned}
\]
and then iterate through other possible combinations of non-zero entries. However, using the importance score matrix $\mathbf{I}$, the least significant entries are identified as $\gamma_{1,1}, \gamma_{2,1},$ and $\gamma_{3,2}$. This allows us to test the reduced problem:

\[
\begin{aligned}
\underset{\gamma}{\min} \quad & \langle \mathbf{C}, \mathbf{\gamma} \rangle, \\
\text{subject to} \quad & \gamma_{1,1} = \gamma_{2,1} = \gamma_{3,2} = 0, \\
                         & \mathbf{\gamma}\mathbf{1} = \mathbf{a}, \quad \mathbf{\gamma}^T\mathbf{1} = \mathbf{b}.
\end{aligned}
\]

If the problem is infeasible, other non-significant entries are selected and tested until a feasible solution is found. The intuition is that with a well-defined importance score reflecting non-zero entries, the algorithm can efficiently move toward the optimal transport plan.

To evaluate this heuristic, we must address three key questions:

\begin{enumerate}
    \item How can we define an effective importance score?
    \item Does this heuristic ensure optimality? In other words, is the cost given by the heuristic algorithm the same as the minimum cost obtained by exhausting all the combinations of the non-zero entries?
    \item Does it improve runtime in practice?
\end{enumerate}

We will explore these questions in the following sections.

\subsubsection{Defining an Importance Score}
This definition of the importance score is not unique, but any valid importance score should indicate the likelihood or necessity of transporting mass through a given entry.

A natural way to define the importance score is as follows:  
\[
\mathbf{I}_{j,k} = \frac{1}{((\mathbf{a} \times \mathbf{b}) \circ \mathbf{C})_{j,k}},
\]
where \( \times \) denotes the outer product and \( \circ \) represents the Hadamard (element-wise) product between matrices.  

The term \((\mathbf{a} \times \mathbf{b})\) reflects the expected transport plan between \(\mu\) and \(\nu\) under a uniform cost assumption, making it a product measure. By applying the Hadamard product with the cost matrix \(\mathbf{C}\), the score is adjusted to reflect both the magnitude of expected mass transported and the associated transportation cost. Intuitively:
\begin{itemize}
    \item Entries expected to transport larger masses are deemed more important.
    \item Cheaper entries (lower \(\mathbf{C}_{j,k}\)) are prioritized.
\end{itemize}

\subsubsection{Optimality of the Heuristic}

The heuristic guarantees optimality only if the importance score perfectly identifies the locations of non-zero entries in the optimal transport plan. For example, if the importance score were ideal, setting non-zero and zero entries accordingly would yield the optimal solution in a single linear programming run.  

In practice, however, the importance score is merely a surrogate for this ideal indicator. As a result, the heuristic can produce suboptimal solutions. Experimental results show that the transport cost with this heuristic can increase by 24\% to 31\% compared to the combinatorial optimum. Despite this, the heuristic consistently outperforms a random selection strategy.

To illustrate: On average, the heuristic yields solutions better than 50\% of all possible transport plans satisfying the sparsity constraints. This means that even with random selections of the non-zero entries guaranteed with a solution, 50\% of the time, the heuristic still produces lower-cost solutions.

Further experimental details can be found in the corresponding section.

\subsubsection{Runtime Efficiency}

The heuristic significantly reduces the runtime of solving sparsity-constrained optimal transport problems, albeit at the cost of suboptimality. For example, with \(n = m = 4\), the heuristic finds a suboptimal solution in approximately 10\% of the time required to exhaustively compute the optimal solution.  

The efficiency gains grow with problem size. Because the traditional solution is combinatorial, as \(n\) and \(m\) increase, the time savings become even more pronounced. Thus, this approach provides a practical trade-off between runtime and solution quality, making it viable for large-scale problems.

\subsection{Experiment Results}
The above discussion is based on the experiment conducted in this Notebook\footnote{\href{https://colab.research.google.com/drive/1N3sWAyYJmGnIINNvNq4yg2stCUrh3z_T?usp=sharing}{Colab Notebook Link}}. In this experiment, we tested our algorithm with $n=m=4$ and 100 randomized sources, sinks, and cost matrices. The experimental results are shown below.

\vspace{0.4cm}
\begin{enumerate}
    \item Surrogate 1:\[\mathbf{I}_{j,k} = ((\mathbf{a} \times \mathbf{b}) \circ \mathbf{C})_{j,k}\]
    \item Surrogate 2:\[\mathbf{I}_{j,k} = \frac{1}{((\mathbf{a} \times \mathbf{b}) \circ \mathbf{C})_{j,k}}\]
    \item Surrogate 3:
\[
\mathbf{I}_{j,k} = \underset{\gamma}{\arg\min} \langle \mathbf{C}, \mathbf{\gamma} \rangle
\quad \text{subject to} \quad \mathbf{\gamma} \mathbf{1} = \mathbf{a}, \quad \mathbf{\gamma}^T \mathbf{1} = \mathbf{b}.
\]
    \item Surrogate 4:\[\mathbf{I}_{j,k} = (\lambda (\mathbf{a} \times \mathbf{b}) + (1-\lambda) \mathbf{C})_{j,k}\] \\ where $\lambda$ = 0.7 in this test.
\end{enumerate}

\begin{table}[H]
\centering
\title{\textbf{Table 2: Experiment Results}}
\begin{tabular}{|c|c|c|c|c|}
\hline
\textbf{Surrogates}          & \textbf{Surrogate 1} & \textbf{Surrogate 2} & \textbf{Surrogate 3} & \textbf{Surrogate 4} \\ 
\hline
{\% of Additional Cost}        & 26\%      & 32\%      & 30\%      & 25\%             \\ 
\hline
{\% of Time Saved}           & 90\%      & 87\%      & 90\%      & 89\%            \\ 
\hline
{\% of Solutions Beat}       & 49\%      & 41\%      & 38\%      & 45\%            \\ 
\hline
\end{tabular}
\caption{The percentages in this table are rounded to the nearest integer.}
\end{table}
Here:
\begin{itemize}
    \item \textbf{\% of Additional Cost}: The percentage increase in transport cost when using the surrogate compared to the optimal solution.
    \item \textbf{\% of Time Saved}: The percentage of time saved by using the surrogate over traditional combinatorial methods.
    \item \textbf{\% of Solutions Beat}: The percentage of solutions from the surrogate that outperforms random selection.
\end{itemize}

\section{Discussion}
Having developed algorithms for both capacity-constrained and sparsity-constrained problems, we now combine them to solve a general real-world problem, as formulated in \eqref{Merged}. First, we apply the algorithms in Section \ref{TOTCC} to compute the transport plan for each time step without the sparsity constraints. Using these plans, we calculate the measures $\mu^i$ and $\nu^i$ for each time step. With the marginals known, we then apply the heuristic algorithm to obtain a suboptimal solution that satisfies the sparsity constraints for each time step. If the heuristic successfully solves the sparsity-constrained problem at every step, we will have a solution to the overall problem \eqref{Merged}.

\section{Conclusion}
In this paper, we extend the classical optimal transport framework by incorporating temporal dynamics to address real-world constraints, such as capacity and sparsity. We propose efficient algorithms for these constraints, achieving significant computational improvements while ensuring practical applicability. Our experiments validate the effectiveness of these methods, demonstrating their potential in dynamic transport scenarios. Future work could focus on integrating additional constraints, enhancing scalability, and refining heuristics to improve accuracy.

\newpage

\bibliographystyle{plain}
\bibliography{citations} 

\end{document}